\documentclass{amsart}
\usepackage{amssymb, amsmath}

\begin{document}
\title{}
\vbox{\hfil {\Large\bf ON EXTENDIBILITY AND DECOMPOSABILITY}\hfil}
\vbox{\hfil {\Large\bf  OF CERTAIN $*$-LINEAR MAPS INTO $\, C ( X )\, $}\hfil}

\author{U. Haag}

\date{\today \\ \texttt{\hfil Contact:haag@mathematik.hu-berlin.de}}
\maketitle
\begin{abstract}
We consider $*$-linear maps into a commutative $C^*$-algebra $C ( X )\, $ of continuous functions on a locally compact Hausdorff space $\, X\, $ with certain specified properties and prove two results: 
\par\noindent
(1) an extension result for a class of $*$-linear maps $\, \varphi : \mathfrak Y \rightarrow C ( X )\, $ which may be called of locally compact type (locally finite in the terminology below) with respect to an inclusion $\, \mathfrak Y \subseteq \mathfrak X\, $ of normed vector spaces, and 
\par\noindent
(2) a minimal decomposition for certain $*$-linear maps into $\, C ( X )\, $ (absolutely continuous maps in our terminology) as a difference of two positive maps.
\end{abstract}
\par\bigskip\noindent
The following Lemma gives a Hahn-Banach-type extension result for certain maps into a commutative 
$C^*$-algebra related to a well known result of Lindenstrauss (cf. \cite{La}, chap. 21, Theorem 1 (10)).  For the relation with the classical theory we also mention the standard reference \cite{Ped}.  Let 
$\, B\, $ denote a commutative (complex) separable $C^*$-algebra and $\, B^{sa}\, ,\, B^+\, $ its selfadjoint (resp. positive) part, $\,\widetilde B\, $ its unitization and $\, C_B\, $ the cone of positive lower semicontinuous functions on the spectrum of $\, B\, $ (in general a function on $\, Spec\, B\, $, continuous or semicontinuous, is meant to be arbitrarily small outside some compact subset). 
Suppose given a separable normed $A^{sa}$-module with $\, A\subseteq B\, $ a suitable $C^*$-subalgebra.  Let $\, C\, $ be a subcone of the cone of positive elements of some (other) normed and ordered $A^{sa}$-module, i.e. it is closed under convex combinations and multiplication by elements from $\, A^+\, $ and contains the order unit. The norm and order of such an $A^{sa}$-module are supposed to be compatible with the norm and order on $\, A\, $ induced by its $C^*$-structure in the obvious way. If $\,\mathfrak X\, $ is an $A^{sa}$-module a $C$-valued (and $A^+$-homogenous) seminorm is a sublinear map 
$\, m : {\mathfrak X} \rightarrow C\, $ which is supposed to satisfy
$$ (i)\quad m ( x + y )\> \leq\> m ( x ) + m ( y )\qquad\qquad\qquad\quad $$
$$ (ii)\quad m ( a\cdot x )\> =\> \vert a\vert \cdot m ( x )\, ,\quad a \in A^{sa}\qquad\qquad\>  . $$
In the following we specialize to the case where $\, A^{sa} = \mathbb R\, $ and $\, C = C_B\, $ consists of positive lower semicontinuous functions on the spectrum of $\, B\, $. Then $\, m\, $ is said to be {\it locally finite} with respect to a closed subspace $\,\mathfrak Y \subseteq \mathfrak X\, $ iff $\, m\, $ is bounded in norm and 
there exists a complementary subspace $\, {\mathfrak Y}^c\, $ with 
$\, \overline{{\mathfrak Y} + {\mathfrak Y}^c} = \mathfrak X\, $,  such that given $\, t_0 \in Spec ( B )\, $ one of the following two conditions applies
\par\medskip\noindent
(i)\quad for every $\, x\in {\mathfrak Y}^c\, $ and $\,\epsilon > 0\, $ there exists a neighbourhood 
$\, {\mathcal U}_{x , \epsilon } ( t_0 )\, $ such that $\, m ( x ) ( t ) < \epsilon\, $ for 
$\, t\in {\mathcal U}_{x , \epsilon } ( t_0 )\, $, in particular $\, m ( x )  ( t_0 ) = 0\, $ for all 
$\, x\in {\mathfrak Y}^c\, $, 
\par\medskip\noindent
(ii)\quad there exists a finitedimensional subspace $\, \mathfrak F \subseteq \mathfrak Y\, $, a complementary subspace  $\, \mathfrak L\, $ with $\, \mathfrak F + \mathfrak L = \mathfrak Y\, $, and a neighbourhood $\, {\mathcal U} ( t_0 )\, $ of $\, t_0\, $ such that for all $\, y \in \mathfrak L\, $ 
$\, m ( y ) ( t ) = 0 \, $ holds for $\, t \in {\mathcal U} ( t_0 )\, $. 
\par\smallskip\noindent
If $\, \phi : \mathfrak Y \rightarrow B^{sa}\, $ is a linear map and $\, m\, $ is a $C_B$-valued seminorm on 
$\,\mathfrak X\, $ dominating $\,\phi\, $ then
$\, m\, $ is said to be {\it balanced} for $\, \phi\, $ if it dominates another seminorm $\, m'\, $ which is locally finite with respect to $\, \mathfrak Y\, $ and dominates $\,\phi\, $, and {\it locally balanced} 
(for $\,\phi\, $) if there exists for every 
$\, t\in Spec\, B\, $ a neighbourhood $\, \mathcal V ( t )\, $ such that 
$\, m {\vert }_{\mathcal V ( t )}\, $ is balanced for $\, {\phi } {\vert }_{\mathcal V ( t )}\, $ (after evaluation of functions on $\, \mathcal V ( t )\, $). 
\par\bigskip\noindent
{\bf Lemma 1.} Assume given a bounded $B^{sa}$-valued and $\mathbb R$-linear map  
$\,\phi\, $ defined on a subspace $\, {\mathfrak Y} \subseteq \mathfrak X\, $ with $\,\mathfrak X\, $ separable, together with a $C_B$-valued 
and ${\mathbb R}^+$-homogenous seminorm $\, m\, $ on $\,\mathfrak X\, $ such that $\, m\, $ is locally balanced for $\, \phi\, $. If 
$\,\delta > 0\, $ is any positive number there exists an $\mathbb R$-linear extension $\,\widetilde\phi\, $ of $\,\phi\, $ to all of 
$\,\mathfrak X\, $, such that $\, \widetilde\phi ( z )\, $ is dominated by $\,  m ( z ) + 
\delta \Vert z \Vert\, $, i.e. 
$\, \vert \widetilde \phi ( x ) ( t ) \vert \leq m ( x ) ( t ) + \delta \Vert x \Vert\, $ for every $\, t\in Spec\, B\, $.
\par\medskip\noindent
{\it Proof.}\quad  Without loss of generality we can assume that $\, m\, $ is locally finite with respect to 
$\, \mathfrak Y\, $ by choosing in case of locally balanced 
$\, m\, $ a locally finite subset $\, \{ \mathcal V ( t_j ) \} \, $ covering $\, Spec\, B\, $ and a corresponding continuous partition of unity $\, \{\, {\lambda }_j\, \}\, $ with 
$\, 0Ê\leq {\lambda }_j ( t ) \leq 1\, ,\, \sum_j\, {\lambda }_j ( t ) = 1\, $ and 
$\, {\lambda }_j ( t ) = 0\, $ for $\, t\notin {\mathcal V} ( t_j )\, $ replacing $\, m\, $ and $\, \phi\, $ by 
$\, ( {\lambda }_j\, m )\, $ and $\, ( {\lambda }_j\, \phi )\, $ and adding up the corresponding extensions.
Let $\, \mathfrak X = \overline{ \mathfrak Y + {\mathfrak Y}^c}\, $ be a corresponding decomposition and 
for any $\, z\in \mathfrak X\, $ and given $\, t\in Spec\, B\, $ let $\, {\overline z}_t\, $ denote its image in the quotient space on dividing
$\, \mathfrak X\, $ by the closure of $\, {\mathfrak Y}^c + {\mathfrak N}_t\, $ with $\, {\mathfrak N}_t\, $ the linear subspace of elements such that $\, m ( x ) ( t ) = 0\, $ for $\, x\in {\mathfrak N}_t\, $. Given 
$\, \delta > 0\, $ define $\, m_{\delta }\, $  by 
$$ m_{\delta } ( z ) ( t )\> =\> m ( z ) ( t )\> +\> \delta\, \Vert {\overline z}_t \Vert \> . $$
and check that this defines a $C_B$-valued seminorm which again is locally finite for the inclusion 
$\, \mathfrak Y \subseteq \mathfrak X\, $. 
Let $\, x \in {\mathfrak Y}^c\, $ be a fixed element and 
$\,\widetilde\phi ( y + a x ) ( t )\, =\, \phi ( y ) ( t ) + a\cdot f ( t )\, $ with 
$\, a\in\mathbb R , f \in B^{sa}\, $ denote an extension of $\,\phi\, $ to the space 
$\, \mathfrak Y + \mathbb R x\, $.
The domination property with respect to $\, m_{\delta }\, $ can be checked pointwise for each 
$\, t \in Spec\, B\, $ and by positive homogeneity of 
$\, m\, $ and $\,\phi\, $ it suffices to check that for each such point 
$$ u_{x , y} := [\phi ( y ) - m_{\delta } ( y - x ) ] ( t )\> <\> f ( t )\> <\> [ - \phi ( z ) + m_{\delta } ( z + x ) ] ( t ) =: l_{x , z} $$
\par\medskip\noindent
holds for all $\, y , z \in \mathfrak Y\, $. The righthandside of this inequality is a lower semicontinuous function which is clearly larger than the upper semicontinuous expression on the left side independent of $\, y\, $ and $\, z\, $, so we may pick $\, f ( t )\, $ from the nonempty closed interval spanned by 
$$ u_x ( t ) := \sup_{y \in \mathfrak Y}\> \{ [ \phi ( y ) - m_{\delta } ( y - x ) ] ( t ) \}\> ,\quad 
l_x ( t ) := \inf_{y \inÊ\mathfrak Y}\> \{ [ - \phi ( y ) + m_{\delta } ( y + x ) ] ( t ) \} . $$
Both expressions are attained inside some ball of finite radius $\, R > 0\, $ around the origin (which can be chosen independent of $\, t\, $) by the relation 
$$ [ - \phi ( r y ) + m_{\delta } ( r y + x ) ] ( t ) \geq [ - \phi ( r y ) + m_{\delta } ( r y ) - m_{\delta } ( - x )] ( t ) $$
\par\smallskip\noindent
where the right side converges to $\, +\infty\, $ for $\, r\to\infty\, $ by positive homogeneity of 
$\, m_{\delta }\, $ and $\,\phi\, $ (unless $\, m_{\delta } ( y ) ( t ) = 0\, $ which implies 
$\, m_{\delta } ( r y + x ) ( t ) = m_{\delta } ( x )\, $). Thus one effectively only needs to consider elements in the compact set $\, {\mathfrak F}^R\, $ which is the intersection of the bounded region as above with the finitedimensional subspace $\,\mathfrak F\, $.
One then shows that by assumption on $\, m\, $ the lefthandside remains upper semicontinuous even after passing to the supremum over all $\, y\in {\mathfrak Y}^R\, $, and that the right side remains lower semicontinuous in the corresponding infimum. 
One notes the following:
if $\, m_{\delta } ( z + x ) ( t ) \geq ( 1 - \epsilon )\, m_{\delta } ( z ) ( t )\, $ the right side of the above inequality is greater or equal to zero, whereas on the corresponding set of points where 
$\, m_{\delta } ( y - x ) ( t ) \geq ( 1 - \epsilon )\, m_{\delta } ( y ) ( t )\, $ the left side of the inequality is less or equal to zero, 
provided that 
$\, ( 1 - \epsilon )\, m_{\delta } ( y ) ( t ) > \phi ( y ) ( t )\, $ is strictly larger.  
Assume that the local infimum 
$\, l_x ( t_0 )\, =\, \inf_{z\in\mathfrak Y}\> l_{x , z} ( t_0 )\, $ is negative at a given point. 
This cannot happen for points with $\, m\, $ subject to condition (i) because $\, m ( x ) ( t_0 ) = 0\, $ implies $\, m ( y + x ) ( t_0 ) = m ( y ) ( t_0 )\, $ for each $\, y \in \mathfrak Y\, $. Then choose a neighbourhood $\, {\mathcal U}\, $ of $\, t_0\, $ such that $\, m ( t )\, $ satisfies condition (ii) in that neighbourhood so the local infimum $\, l_x ( t )\, $ is attained within the compact set 
$\, {\mathfrak F}^R\, $. Assume that $\, l^R_x\, $ is not lower semicontinuous at 
$\, t_0\, $. Then there exists $\, c > 0\, $ such that in every neighbourhood of $\, t_0\, $ there are points satisfying $\, l_x ( t_0 ) - l_x ( t ) > c\, $. 
Choose a sequence of points 
$\, \{ t_n \}\, $ converging to $\, t_0\, $ such that $\, l_x ( t_0 ) - l_x ( t_n ) > c\, $ for every $\, n\, $. Then there is a sequence 
$\,\{ z_n \}\, $ inside the compact region as above such that 
$\,   l_{x , z_n} ( t_n ) - l_x ( t_n ) \leq 1 / n\, $ for all $\, n\, $.
The sequence $\, \{ z_n \}\, $ has an adherence point $\, z_\infty\, $ and since this implies that 
$\, l_{x , z_\infty } ( t )\, $ is an adherence point for $\, l_{x , z_n} ( t )\, $ at each point 
$\, t\, $  one concludes that for every fixed $\, t\, $ there exists a certain subsequence of 
$\,\{ l_{x , z_n} ( t ) \}\, $
which converges to $\, l_{x , z_\infty } ( t )\, $. By a diagonal method one constructs a subsequence such that 
$\, \vert l_{x , z_n} ( t_n ) - l_{x , z_\infty } ( t_n ) \vert < 1 / n\, $. Then 
$\, z_\infty\, $ satisfies 
$\, l_{x , z_\infty } ( t_0 ) \geq l_x ( t_0 )\, $ and hence
$\, l_{x , z_\infty } ( t_0 ) - l_{x , z_\infty } ( t_n ) > c - 1 / n\, $ which is impossible since it would make 
$\,  l_{x , z_\infty }\, $ lower semidiscontinuous at $\, t_0\, $. 
Thus for all values of $\, t\, $ where $\, m\, $ is subject to condition (ii) as above, in particular at all places where $\, l_x ( t )\, $ is negative, the function 
$\, l_x\, $ is lower semicontinuous at $\, t\, $. In the same way one shows that the function $\, u_x\, $ is upper semicontinuous at all points $\, t\, $ where $\, u_x ( t ) > 0\, $ (resp. points where $\, m\, $ satisfies condition (ii)). If on the other hand $\, m\, $ is satisfies condition (i) at a given $\, t\, $ one gets 
$\, u_x ( t ) = 0 = l_x ( t )\, $. Since $\, m ( x )\, $ is continuous at $\, t_0\, $ the value of $\, m ( x ) ( t )\, $ becomes arbitrarily small in some neighbourhood of 
$\, t_0\, $, making the expressions $\, u_x\, $ and $\, l_x\, $ continuous at $\, t_0\, $. It follows from the Continuous Selection Criterion of \cite{Mi} that one may choose a continuous function 
$\, f \, $ satisfying 
$$ u_x ( t ) \> \leq\> f ( t )\> \leq\> l_x ( t ) \> . $$
Then the assignment 
$$ {\widetilde\phi} ( y + a x ) ( t ) \> =\> \phi ( y ) ( t ) + a f ( t )  $$
gives an extension of $\, \phi\, $ to the subspace $\, \mathfrak Y + \mathbb R x\, $ dominated by 
$\, m_{\delta }\, $.  One proceeds by induction. Since $\, \delta\, $ can be chosen arbitrarily small in the above argument and $\, m_{\delta }\, $ is locally finite for the inclusion $\, \mathfrak Y + \mathbb R x \subseteq \mathfrak X\, $ by choosing a suitable complement of $\, \mathbb R x\, $ in 
$\, {\mathfrak Y}^c\, $ one can repeat the argument extending $\, \widetilde\phi\, $ to a larger subspace 
$\, \mathfrak Y + \mathbb R x + \mathbb R x'\, $ with $\, x , x'\in {\mathfrak Y}^c\, $ spanning a twodimensional subspace, but this time with the extension dominated by $\, m_{\delta } ( z ) ( t ) + 
{\delta } / 2 \Vert {\widetilde z}_t \Vert\, $ where $\, {\widetilde z}_t\, $ is the image of $\, z\, $ in the quotient of $\,\mathfrak X\, $ divided by the closure of the complementary subspace 
$\, {{\mathfrak Y}'}^c\, $ of $\, \mathfrak Y + \mathbb R x\, $ plus the nilspace of $\, m ( t )\, $. Inductively one extends $\,\phi\, $ to a linear map 
$\,\widetilde\phi\, $ dominated by $\, m + 2 \delta \Vert\cdot\Vert \, $ on a dense linear subspace of   
$\,\mathfrak X\, $, by separability, and by continuity to all of $\,\mathfrak X\, $  \qed
\par\bigskip\noindent
{\it Example.}\quad Assume given a linear map $\,\phi\, $ into $\, B^{sa}\, $ defined on a real subspace 
$\,\mathfrak Y \subseteq \mathfrak X\, $ and a $C_B$-valued seminorm $\, m\, $ on $\,\mathfrak X\, $ dominating $\,\phi\, $ such that $\, m\, $ satisfies the second condition in the definition of a locally finite seminorm for the subspace $\,\mathfrak Y\, $ in all but one point $\, t_0\, $. From the ordinary 
Hahn-Banach Theorem one gets an extension $\, {\widetilde\phi }_{t_0} ( x )\, $ of 
$\, \phi ( x ) ( t_0 )\, $ to all of $\,\mathfrak X\, $ dominated by $\, m ( x ) ( t_0 )\, $. The kernel of 
$\, {\widetilde\phi }_{t_0}\, $ decomposes into its intersection with $\,\mathfrak Y\, $ and a complementary part which is denoted $\, {\mathfrak Y}^c\, $. Choose an increasing sequence of finitedimensional subspaces $\, {\mathfrak F}_n\subseteq \mathfrak Y\, $ with dense union in $\,\mathfrak Y\, $ and for 
given $\, \delta > 0\, $ consider the quotient seminorms 
$$ {\overline m}_{\delta } ( x ) ( t )\> =\> \inf_{y\in {\mathfrak Y}^c + {\mathfrak N}_t}\,\left[  m ( x + y ) ( t )\> +\> \delta \Vert x + y \Vert \right] \> , $$
$$ {\widetilde m}_{\delta } ( x ) ( t )\> =\> m ( x ) ( t )\, +\, \inf_{y\in {\mathfrak Y}^c + {\mathfrak N}_t}\, 
\left[ m ( y ) ( t )\, +\, \delta \Vert x + y \Vert \right] \> . $$
Since $\, {\widetilde\phi }_{t_0}\, $ factors over the quotient by $\, {\mathfrak Y}^c\, $ one gets that 
$\, \phi ( x ) ( t_0 )\, $ is (strictly) dominated by $\, {\overline m}_{\delta } ( x ) ( t_0 )\, $ for every 
$\, x\in\mathfrak Y\, $. Then for given $\, {\mathfrak F}_n\, $ there exists a neighbourhood 
$\, {\mathcal U}_n ( t_0 )\, $ of $\, t_0\, $ such that $\, \phi ( x ) ( t )\, $ is still (strictly) dominated by 
$\, {\overline m}_{\delta } ( x ) ( t )\, $ for all $\, x\in {\mathfrak F}_n\, $. Define a new $C_B$-valued seminorm on $\, {\mathfrak F}_n + {\mathfrak Y}^c\, $ by 
$$ m_{\delta , n} ( x ) ( t )\> =\quad {\overline m}_{\delta } ( x ) ( t )\qquad  if\quad
t\in \overline{\mathcal U}_n ( t_0 ) \> , $$
$$ m_{\delta , n} ( x ) ( t )\> =\quad {\widetilde m}_{\delta } ( x ) ( t )\qquad otherwise
\qquad\>  .  $$ 
Then inductively define a $C_B$-valued seminorm $\, m_{\delta }'\, $ in the following manner: suppose that $\, m_{\delta }'\, $ has already been constructed for all $\, y\in {\mathfrak F}_n + {\mathfrak Y}^c\, $, and for $\, x\in {\mathfrak F}_{n+1} + {\mathfrak Y}^c\, $ put
$$ m_{\delta }' ( x ) ( t )\> =\> \inf_{y\inÊ{\mathfrak F}_n}\, \left[ m_{\delta }' ( y ) ( t ) + 
m_{\delta , n+1} ( x - y ) ( t )Ê\right] $$
and check that this defines a locally finite seminorm for the inclusion 
$\,\mathfrak Y\subseteq \mathfrak X\, $ dominating $\,\phi\, $ and dominated by 
$\, {\widetilde m}_{\delta }\, $. Thus up to an arbitrarily small correction the original seminorm $\, m\, $ is balanced for $\, \phi\, $. 
Then by the same argument if $\, m\, $ is subject to condition $\, (ii)\, $ of the above definition in all but finitely many places (resp. on a discrete nowhere dense subset) of $\, Spec \, B\, $ then up to an infinitesimal correction term $\, m\, $ is locally balanced. This property will also be called 
{\it locally $\delta $-balanced}.
\par\bigskip\noindent
Let $\, A\, $ be a $C^*$-algebra. A linear (resp. sublinear) map 
$\,\phi : A^{sa} \rightarrow B^{sa}\, $ is said to be {\it absolutely continuous} iff it is normbounded and if the pointwise norm 
$$ \Vert \phi \Vert ( t ) := \sup_{\Vert x \Vert \leq 1}\> \vert \phi ( x ) ( t ) \vert $$
with $\, t\in Spec\, B\, $ is an element of $\, \widetilde B^+\, $ (i.e. continuous and converging at infinity). 
 \par\bigskip\noindent
{\bf Lemma 2.} (Jordan decomposition) 
Let $\, A\, $ be a unital $C^*$-algebra and $\, B\, $ commutative. Any selfadjoint absolutely continuous linear map $\, \phi :  A \rightarrow B\, $ can be uniquely written as a difference 
$\, \phi = {\phi }_+ - {\phi }_-\, $ of two positive (absolutely continuous) linear maps 
$\,{\phi }_+\, ,\, {\phi }_-\, :\, A \rightarrow  \widetilde B\, $ such that 
$\,\Vert \phi \Vert\, ( t )\> =\> 
\Vert {\phi }_+ \Vert\, ( t ) + \Vert {\phi }_- \Vert\, ( t )\, $ holds for each $\, t\in Spec\, B \, $. 
\par\medskip\noindent
{\it Proof.}\quad  Consider the $w^*$-continuous family 
$\,\{ {\phi }_t\,\vert\, t \in Spec\, B  \}\, $ of selfadjoint linear functionals given by composing $\,\phi\, $ with the evaluation map at $\, t \in Spec\, B\, $, and the two corresponding families 
$\,\{ {\phi }_{t , +} \}\, ,\, \{ {\phi }_{t , -} \}\, $ of positive functionals such that 
$$ {\phi }_t ( x )\> =\> {\phi }_{t , +} ( x )\> -\> {\phi }_{t , - } ( x ) $$
\par\smallskip\noindent
equals the Jordan decomposition of $\, {\phi }_t\, $. We claim that these families are again 
$w^*$-continuous, and hence glue together to positive linear maps 
$$ {\phi }_+\, :\, A\>\rightarrow\> B\qquad ,\qquad 
      {\phi }_-\, :\, A\>\rightarrow\> B $$
\par\smallskip\noindent
such that $\, \Vert \phi \Vert ( t )\, =\, \Vert {\phi }_+ \Vert ( t )\, +\, \Vert {\phi }_- \Vert ( t )\, $. We first do the case where $\, A\, $ is commutative which is especially simple (and instructive). Since the norm of 
$\, {\phi }_{\pm , t} \, $ is given by $\, {\phi }_{\pm , t} ( 1 )\, $ at each point $\, t\, $, the difference of the normfunctions as well as their sum is continuous, because 
$$ \Vert {\phi }_{+ , t } \Vert \, -\, \Vert {\phi }_{- , t} \Vert\> =\> {\phi } ( 1 ) ( t )\> ,\quad 
\Vert {\phi }_{+ , t}Ê\Vert\, +\, \Vert {\phi }_{- , t} \Vert\> =\> \Vert {\phi }_t \Vert  $$
so that the pointwise Jordan decomposition applied to the unit element is continuous. Let 
$\, 0 \leq h\leq 1\, $ be an arbitrary positive element smaller than $\, 1\, $. Consider the linear map 
$\, {\psi }_h : A \rightarrow B\, $ given by $\, {\psi }_h ( x ) = \phi ( h x )\, $. We claim that 
$\, {\psi }_{h , \pm , t} ( x ) = {\phi }_{\pm , t} ( h x )\, $. It is clear that the right side of this equation gives a  positive decomposition of $\, {\psi }_{h , t}\, $ so that by the special properties of the Jordan decomposition one gets 
$$ {\phi }_{\pm , t} ( h )\> =\> \Vert {\phi }_{\pm , t} ( h\cdot )\Vert\> 
\geq\> \Vert {\psi }_{h , \pm , t} \Vert \> . $$
Given $\, \epsilon > 0\, $ and $\, t\in Spec\, B\, $ there exists a positive element 
$\,  K_{\epsilon , t}\in A\, $ of norm less or equal to one such that 
$$  {\phi }_ {+ , t} ( x K_{\epsilon , t} )   \leq \epsilon \Vert x\Vert \> ,
\quad {\phi }_{- , t} ( x ( K_{\epsilon , t} - 1 ) )\> \leq\> \epsilon \Vert x\Vert \> . $$ 
Then 
$$ {\phi }_{+ , t} ( h )\> +\> {\phi }_{- , t} ( h )\> \leq\> {\phi }_t ( h ( 1 - K_{\epsilon , t} ) )\> -\> 
{\phi }_t ( h K_{\epsilon , t} )\> +\> 2\epsilon\> =\> {\psi }_{h , t} ( 1 - 2 K_{\epsilon , t} ) \> +\> 2\epsilon\>\leq\> \Vert {\psi }_{h , t} \Vert\> +\> 2\epsilon $$ 
so that the reverse inequality 
$$   {\phi }_{\pm , t} ( h )\>\leq\> \Vert {\psi }_{h , \pm , t} \Vert \>  $$
follows in the limit $\, \epsilon \to 0\, $. In particular both families 
$\, \{ {\phi }_{\pm , t} ( h ) \}\, $ are lower semicontinuous, since their difference is continuous, and their sum is lower semicontinuous being the norm function of the continuous linear map $\, {\psi }_h\, $. The same holds with $\, h\, $ replaced by $\, 1 - h\, $. Then $\, {\phi }_{\pm , t} ( h )\, $ and 
$\, {\phi }_{\pm , t} ( 1 - h )\, $ are both lower semicontinuous functions adding up to the 
continuous function $\, {\phi }_{\pm , t} ( 1 )\, $ which means that both functions must be continuous by themselves. Since any element is a linear combination of positive elements one concludes that 
the families $\, \{ {\phi }_{\pm , t} ( x ) \}\, $ are continuous for any $\, x\in A^{sa}\, $.
Note that weakening the assumption to $\,\phi\, $ being only absolutely $\delta $-continuous, the same argument will produce a pointwise 
Jordan decomposition which is $\delta $-continuous in each single element, i.e. 
$$ \limsup_{t\to t_0}\, \vert {\phi }_{\pm , t} ( x ) - {\phi }_{\pm , t_0} ( x ) \vert\> \leq\> \delta \Vert x\Vert $$
The case of noncommutative $\, A\, $ is a bit more involved, although the argument in case of the unit element is just the same, the problem being that the map $\, x\mapsto h x\, $ with $\, 0 \leq h < 1\, $ is no longer positive in the noncommutative case.  We will use an approximation by finitedimensional real operator subsystems of $\, {\mathfrak F}_n\subseteq A^{sa}\, $ such that the union of the 
$\, \{ {\mathfrak F}_n \}\, $ is dense in $\, A^{sa}\, $. 
Let $\, {\mathfrak F}_c \subseteq A\, $ be any (complex) finitedimensional operator subsystem with selfadjoint part $\, \mathfrak F\, $ and $\, {\mathcal P}_{\mathfrak F} \subseteq 
 ({\mathfrak F}_c)^*\, $ the closure of the set of pure or extremal states, consisting of extremal positive functionals of norm one, whereas the (closed) subset of arbitrary positive functionals of norm one is denoted by  $\, {\mathcal S}_{\mathfrak F}\, $. The Kadison function representation (cf. \cite{Kad})represents $\, \mathfrak F\, $ isometrically as the subspace $\, {\mathit A} ( \mathfrak F )\, $ of continuous realvalued functions on $\, {\mathcal P}_{\mathfrak F}\, $ which admit an affine extension to 
 $\, {\mathcal S}_{\mathfrak F}\, $, therefore also contained in 
 $\, C ( {\mathcal S}_{\mathfrak F} )\, $ which in turn is naturally contained in $\, C ( {\mathcal S}_A )\, $, the commutative $C^*$-algebra of continuous functions on the state space of $\, A\, $. 
Given $\, \mathfrak F\, $ as above
one may extend the restriction $\, \phi {\vert }_{\mathfrak F}\, $ to a linear map 
$$ {\phi }^{\delta }_{\mathfrak F} :\> C ( {\mathcal S}_A )^{sa}\> \longrightarrow\> B^{sa} $$
with 
$\, \Vert {\phi }^{\delta }_{\mathfrak F} \Vert ( t ) \leq \Vert \phi {\vert }_{\mathfrak F} \Vert ( t ) + \delta \, $ for any given $\,\delta > 0\, $ from Lemma 1 (due to finitedimensionality of $\,\mathfrak F\, $) since 
$\,\phi {\vert }_{\mathfrak F}\, $ is dominated by the locally finite and continuous seminorm 
$\, m ( x ) ( t ) = \Vert \phi {\vert }_{\mathfrak F} ( t ) \Vert  \Vert x \Vert\,  $. This may be done in such a way that the norm function of the extension is continuous. For an increasing sequence of finitedimensional operator subsystems $\, \{ {\mathfrak F}_n \}\, $ with dense union in $\, A^{sa}\, $ as above one may choose a corresponding sequence $\,  ( {\delta }_n )_n\, $ converging to zero.
Then since the norm functions 
$\, \{ \Vert {\phi }^{{\delta }_n }_{{\mathfrak F}_n , t} \Vert \}\, $ are continuous and monotonously increasing (up to the value of $\, {\delta }_n\, $) to the continuous limit function  $\, \Vert {\phi }_t \Vert\, $ the sequence converges uniformly. 
For our purposes we need to refine this construction. 
For any fixed $\, t\, $ consider an arbitrary extension 
$$ {\widetilde\phi }_t :\> C ( {\mathcal S}_A )^{sa} \>\longrightarrow\> \mathbb R $$
of $\, {\phi }_t\, $ satisfying $\, \Vert {\widetilde\phi }_t \Vert = \Vert {\phi }_t \Vert\, $ which exists by the Hahn-Banach Theorem and put $\, {\widetilde\psi }_{t , h} ( x ) = {\widetilde\phi }_t ( h x )\, $. Since the norm is the same as the norm of $\,{\phi }_t\, $ the Jordan decomposition of $\, {\widetilde\phi }_t\, $ restricts to the Jordan decomposition of $\, {\phi }_t\, $ on $\, A^{sa}\, $. Therefore by the argument above one gets 
$$ \Vert {\widetilde\psi }_{t , h , \pm } \Vert\> =\> {\widetilde\phi }_{t , \pm } ( h )\> =\> 
{\phi }_{t , \pm } ( h ) \> . $$
For $\, x\in C ( {\mathcal S}_A )^{sa}\, $ define
$$ \overline\phi ( x ) ( t )\> =\> \sup\> \bigl\{\, {\widetilde\phi }_t ( x )\, \bigm\vert\, {\widetilde\phi }_t {\vert }_{A^{sa}} \equiv {\phi }_t\, ,\, \Vert {\widetilde\phi }_t \Vert = \Vert {\phi }_t \Vert \bigr\}  $$
and 
$$ \underline\phi ( x ) ( t )\> =\> -\inf\> \bigl\{\, {\widetilde\phi }_t ( x )\,\bigm\vert\, 
{\widetilde\phi }_t {\vert }_{A^{sa}} \equiv {\phi }_t\, ,\, \Vert 
{\widetilde\phi }_t \Vert = \Vert {\phi }_t \Vert \bigr\} \>  , $$
so that $\, \overline\phi ( x ) ( t )\, $ corresponds to the function $\, l_x ( t )\, $ in the proof of Lemma 1 and $\, \underline\phi ( x ) ( t )\, $ corresponds to the function $\, - u_x ( t )\, $ with $\, m_{\delta } \, $ replaced by the seminorm $\, \Vert {\phi }_t \Vert\, \Vert x \Vert\, $. From this one easily finds that both functions are upper semicontinuous (since $\,\Vert {\phi }_t \Vert\, $ is continuous). Also define 
$$ {\overline\phi }_n ( x ) ( t )\> =\> \sup\, \bigl\{ {\phi }_{{\mathfrak F}_n} ( x ) ( t )\,\bigm\vert\, 
{\phi }_{{\mathfrak F}_n} {\vert }_{{\mathfrak F}_n} \equiv {\phi } {\vert }_{{\mathfrak F}_n}\, ,\, 
\Vert {\phi }_{{\mathfrak F}_n} ( t ) \Vert \leq \Vert {\phi } ( t ) \Vert + {\delta }_n\,\bigr\} \>, $$
$$ {\underline\phi }^n ( x ) ( t )\> =\> - \inf\, \bigl\{\, {\phi }_{{\mathfrak F}_n} ( x ) ( t )\,\bigm\vert\, 
{\phi }_{{\mathfrak F}_n} {\vert }_{{\mathfrak F}_n} \equiv {\phi } {\vert }_{{\mathfrak F}_n}\, ,\, 
\Vert {\phi }_{{\mathfrak F}_n} ( t ) \Vert \leq \Vert {\phi } ( t ) \Vert + {\delta }_n \,\bigr\}Ê\> , $$
i.e. the $*$-linear maps $\, {\phi }_{{\mathfrak F}_n} : C ( {\mathcal S}_A )^{sa} \rightarrow B^{sa}\, $ are assumed to be extensions of the restriction of $\, \phi\, $ to the corresponding finitedimensional operator subsystems $\, \{ {\mathfrak F}_n \}\, $ which are dominated by the slightly larger but continuous seminorm $\, m_n ( x ) ( t ) = \Vert {\phi } ( t ) \Vert\, \Vert x \Vert + {\delta }_n \Vert x \Vert\, $ instead of 
choosing $\, ( \Vert \phi {\vert }_{{\mathfrak F}_n} \Vert + {\delta }_n ) \Vert x \Vert\, $ as suggested above. Since the norm functions $\, \{ \Vert {\phi }_{{\mathfrak F}_n} ( t )  \Vert \}\, $ are uniformly convergent to 
$\, \Vert {\phi } ( t ) \Vert\, $ in the latter case, the same is true for the norm functions of extensions dominated by $\, m_n ( x ) ( t )\, $ which are therefore $\epsilon $-continuous whenever 
$\, n \geq N_{\epsilon }\, $ is large enough with arbitrary $\,\epsilon > 0\, $ (and in any case lower semicontinuous), although they need not be exactly continuous as in the case considered above. This implies that also the Jordan decompositions $\, {\phi }_{{\mathfrak F}_n , \pm } ( x ) ( t )\, $ are 
$\epsilon $-continuous for $\, n \geq N_{\epsilon }\, $ from the proof of the case of abelian $\, A\, $ above.
From the definitions it is clear that the functions $\, \{Ê{\overline\phi }_n ( x ) ( t ) \}\, $ are pointwise monotonously decreasing towards $\, \overline\phi ( x ) ( t )\, $, with 
$\, \{ {\underline\phi }^n ( x ) ( t ) \}\, $ monotonously increasing towards $\, \underline\phi ( x ) ( t )\, $ (as shown in the proof of Lemma 1 the functions 
$\, {\overline\phi }_n ( x ) ( t )\, ,\, {\underline\phi }^n ( x ) ( t )\, $ are continuous and hence arise as actual values of a suitable $\, {\phi }_{{\mathfrak F}_n} ( x ) ( t )\, $). For fixed $\, t\in Spec\, B\, ,\, \epsilon > 0\, $ and any extension $\, {\widetilde\phi }_t\, $ of $\, {\phi }_t\, $ as above there exists
$\, 0\leq K_{t , \epsilon }\leq 1\, $ in $\, C ( {\mathcal S}_A )^{sa}\, $ such that 
$$ {\widetilde\phi }_t ( h K_{t , \epsilon } )\>\leq\> {\phi }_{t , +} ( h )\>\leq\> {\widetilde\phi }_t ( h K_{t , \epsilon } ) + \epsilon\> $$
for each $\, h\geq 0\, $. Then for each $\, n\, $ choose $\, {\phi }_{{\mathfrak F}^t_n}\, $ with 
$\, {\phi }_{{\mathfrak F}^t_n} ( h K_{t , \epsilon } ) ( t ) \geq {\widetilde\phi }_t ( h K_{t , \epsilon } )\, $ implying 
$$ {\phi }_{t , +} ( h ) - \epsilon\> \leq\> {\widetilde\phi }_t ( h K_{t , \epsilon } )\>\leq\> 
{\phi }_{{\mathfrak F}^t_n} ( h K_{t , \epsilon } ) ( t )\>\leq\> 
{\phi }_{{\mathfrak F}^t_{n , +}} ( h K_{t , \epsilon } ) ( t )\>\leq\> {\phi }_{{\mathfrak F}^t_{n , +}} ( h ) ( t ) $$ 
since $\, \Vert K_{t , \epsilon } \Vert \leq 1\, $. For any local supremum $\, t_s\, $ of the function 
$\, {\phi }_+ ( h ) ( t )\, $ (i.e. $\, {\phi }_+ ( h ) ( t_s ) \geq \limsup_{t\to t_s}\, {\phi }_+ ( h ) ( t )\, $)
there is a neighbourhood $\, {\mathcal U}_{t_s}\, $ of the point $\, t_s\, $ such that 
$\, {\phi }_{{\mathfrak F}^{t_s}_{n , +}} ( h ) ( t ) \geq {\phi }_+ ( h ) ( t ) - 2 \epsilon\, $ for $\, t\in {\mathcal U}_{t_s}\, $ by continuity of the function 
$\, {\phi }_{{\mathfrak F}^{t_s}_n} ( h K_{t_s , \epsilon } ) ( t ) \leq 
{\phi }_{{\mathfrak F}^{t_s}_{n , +}} ( h ) ( t )\, $. The union 
$\, \mathcal U = \bigcup\, {\mathcal U}_{t_s}\, $ is dense in $\, X = Spec\, B\, $, and the covering by the open sets 
$\, {\mathcal U}_{t_s}\, $ admits a locally finite subcovering which for simplicity is still denoted 
$\, \{ {\mathcal U}_{t_s} \}\, $. For $\, t\in\mathcal U\, $ define 
$$ {\phi }_{n , t} ( x )\> =\> \min\, \{ {\phi }_{{\mathfrak F}^{t_s}_n} ( x ) ( t )\,\vert\, t\in {\mathcal U}_{t_s} \}
\> . $$
If $\, t_0\notin\mathcal U\, $ consider the net of functionals 
$\, \{ {\phi }_{n , t} ( x ) ( t_0 )\,\vert\, t\to t_0\, ,\, t\in\mathcal U \}\, $ which is is uniformly bounded, hence admits a subnet converging in the $w^*$-topology to a limit functional $\, {\phi }_{n , t_0}\, $ such that 
$\, {\phi }_{n , t_0, +} ( h ) = \limsup_{t\to t_0}\, {\phi }_{n , t , +} ( h ) ( t_0 )\, $ by choosing the convergent subnet subject to the condition 
$\, \Vert {\psi }_{n , h , t_0} \Vert = 
\limsup_{t\to t_0}\, \Vert {\psi }_{n , h , t} \Vert\, $ with 
$\, {\psi }_{n , h , t_0} ( x ) = {\phi }_{n , t_0} ( h x )\, $ etc.. One verifies that by this definition 
$\, {\phi }_{n , t , +} ( h )\, $ is an upper $\epsilon $-semicontinuous function and majorizes 
$\, {\phi }_{t , +} ( h ) - 2 \epsilon\, $. In order that 
$\, \{ {\phi }_{n , t , +} ( h ) {\}}_n\, $ converges pointwise to the value of $\, {\phi }_{t , +} ( h )\, $ it is sufficient that at each point the functionals $\, \{ {\phi }_{n , t} \}\, $ converge to $\, {\phi }_t\, $, and that the pointwise norms $\, \{ \Vert {\phi }_{n , t} \Vert \}\, $ converge to the value of $\, \Vert {\phi }_t \Vert\, $. The first convergence is immediate since for $\, x\in {\mathfrak F}_n\, $ one has 
$\, {\phi }_{n , t} ( x ) = {\phi }_t ( x )\, $. The second property then follows from the global estimate 
$$ \Vert {\phi }_{n , t} \Vert\> \leq\> \Vert {\phi }_t \Vert + {\delta }_n $$
which holds for all $\, n\in\mathbb N\, $ due to the construction of $\, {\phi }_{n , t}\, $. Then 
$\, {\phi  }_{t , +} ( h )\, $ is the decreasing limit of upper $\epsilon $-semicontinuous functions, hence must also be upper $\epsilon $-semicontinuous. As $\, \epsilon > 0\, $ is arbitrary one concludes that 
it is in fact upper semicontinuous. Then the same argument applies to $\, {\phi }_{t , +} ( 1 - h )\, $ and since the sum of these functions is continuous, both must be continuous, hence also 
$\, {\phi }_{t , -} ( h )\, $ and $\, {\phi }_{t , -} ( 1 - h )\, $. Since $\, A\, $ is linearly generated by its positive elements the assignment 
$$ x\>\mapsto\> {\phi }_{t , \pm} ( x )\> ,\quad x\in A  $$
is continuous, i.e. $\, {\phi }_{t , \pm} ( x )\, $ defines an element of $\, B\, $\qed
\par\bigskip\noindent
{\it Remark.}\quad ( i ) If $\, A\, $ is a $B$-algebra, i.e. if $\, B\, $ is a central subalgebra of $\, A\, $ and 
$\,\phi\, $ is an absolutely continuous  $B$-linear map, it is plain to see that $\, {\phi }_{\pm }\, $ are also 
$B$-linear. Namely, the functional $\, {\phi }_t\, $ will factor over the evaluation map 
$\, A \twoheadrightarrow A_t\, $ dividing by the ideal $\, C_0 ( Spec ( B )\backslash \{ t \} )\cdot A\, $, so that its Jordan decomposition also factors over this quotient. This implies that $\, {\phi }_{\pm }\, $ are 
$B$-linear.
\par\smallskip\noindent
( ii ) Since any linear map $\, \phi : A^{sa} \rightarrow B^{sa}\, $ can be elementwise approximated by 
suitable extensions 
$\,Ê {\phi }^{\delta }_{{\mathfrak F}_{\lambda }}\, $ of $\, \phi {\vert }_{{\mathfrak F}_{\lambda }}\, $
with $\, \{ {\mathfrak F}_{\lambda } {\}}_{\lambda\in\Lambda }\, $ a cofinal net of finitedimensional operator subsystems of $\, A^{sa}\, $ as in the proof above, with continuous norm-function 
$\, \Vert {\phi }^{\delta }_{{\mathfrak F}_{\lambda }} \Vert ( t )\, $ (the normfunction of an arbitrary bounded map on a finite dimensional space is continuous) one gets that an arbitrary $*$-linear map into a commutative $C^*$-algebra is the elementwise limit of (Jordan) decomposable maps (compare \cite{Haagerup}).
\par

\end{document}